\documentstyle[12pt,psfig,amssymb]{amsart}

\newcommand{\new}{\newcommand}
\new{\tnsr}{\otimes}
\new{\tensor}{\otimes}
\new{\superset}{\supset}
\new{\union}{\cup}
\new{\CC}{\Bbb C}
\new{\NN}{\Bbb N}
\new{\RR}{\Bbb R}
\new{\ZZ}{\mbox{$\Bbb Z$}}
\new{\FF}{\Bbb F}
\new{\TT}{\Bbb T}
\new{\tr}{\operatorname{tr}}
\new{\abs}[1]{\left|#1 \right|}
\new{\norm}[1]{\left\| #1 \right\|}
\new{\bracket}[1]{\langle #1 \rangle}
\new{\defequals}{\stackrel{\rm def}{=}}
\new{\comb}[2]{{#1 \choose #2}}
\new{\lbl}[1]{\label{#1}
        \if\draft y
                \smash{\makebox[0pt]{\hspace{-0.5in}
                        \raisebox{8pt}{\rm\tiny #1}}}
        \fi
}

\newtheorem{proposition}{Proposition}
\newtheorem{theorem}{Theorem}

\theoremstyle{definition}
\newtheorem{lemma}{Lemma}
\newtheorem{claim}{Claim}
\newtheorem{corollary}{Corollary}
\newtheorem{conjecture}{Conjecture}
\newtheorem{definition}{Definition}

\theoremstyle{remark}

\newcounter{letter}
\newcounter{numeral}
\newcounter{Numeral}

\newenvironment{romanlist}{
\begin{list}{(\roman{numeral})}{\usecounter{numeral}}
}{\end{list}}

\new{\prop}[1]{\begin{proposition} \lbl{#1}}
\new{\thm}[1]{\begin{theorem} \lbl{#1}}
\new{\lem}[1]{\begin{lemma} \lbl{#1}}
\new{\clm}[1]{\begin{claim} \lbl{#1}}
\new{\cor}[1]{\begin{corollary} \lbl{#1}}
\new{\cnj}[1]{\begin{conjecture} \lbl{#1}}
\new{\fig}{\begin{figure}[hbt]}
\new{\eq}[1]{\begin{equation} \lbl{#1}}

\new{\eprop}{\end{proposition}}
\new{\ethm}{\end{theorem}}
\new{\elem}{\end{lemma}}
\new{\eclm}{\end{claim}}
\new{\ecor}{\end{corollary}}
\new{\ecnj}{\end{conjecture}}
\new{\efig}[2]{\caption{#1} \label{#2} \end{figure}}
\new{\eeq}{\end{equation}}

\new{\pic}[5]{\raisebox{#3pt}{
\hspace{#4pt}\psfig{file=#1.ps,height=#2pt,silent=}\hspace{#5pt}}}

\begin{document}
\author{Stephen Sawin}
\address{Math Dept. 2-265, MIT,
Cambridge, MA 02139-4307, (617) 253-4995}
\email{sawin@@math.mit.edu}

\new{\cob}{\operatorname{Cob}}
\new{\hil}{\operatorname{Hil}}
\new{\draft}{n}

\title{Direct Sum Decompositions and \\Indecomposable TQFT's}

\maketitle
\begin{abstract}
The decomposition of an arbitrary axiomatic topological quantum field
theory or TQFT into indecomposable theories is given.  In particular,
unitary TQFT's in arbitrary dimensions are shown to decompose into a
sum of theories in which the Hilbert space of the sphere is
one-dimensional, and indecomposable two-dimensional theories are
classified.
\end{abstract}

\section{Introduction}

In \cite{DJ94}, Durhuus and Jonsson define a notion of direct sum of
axiomatic topological quantum field theories, or TQFT's.  They show
that every unitary TQFT in two dimensions can be written as a direct
sum of theories in which  the Hilbert space associated to the circle is
one-dimensional.  Such theories are easily classified and are
described in terms of Euler number.

Durhuus and Jonsson leave  open two questions:  Higher dimensional
theories and nonunitary theories.  The first they explicitly address,
suggesting every unitary theory in $d$ dimensions can be decomposed
into a direct sum of theories in which the Hilbert space associated to the
sphere $S^{d-1}$ is one-dimensional.  In this paper we give a
complete decomposition theory for TQFT's over an algebraically closed
field, and describe the indecomposable theories explicitly in two
dimensions.  In particular we prove the conjecture suggested by Durhuus and
Jonsson.

The nonunitary, indecomposable TQFT's in two dimensions which we
construct are remarkably degenerate:  Any cobordism of genus two or
higher gets sent to the zero operator on the appropriate space.   In
particular, we construct many counterexamples to the conjecture that a
TQFT is determined by its values on closed manifolds.  The same
conjecture for unitary TQFT's is still open.

This paper is divided into four sections after this introduction.
Section 2 defines TQFT's and the direct sum operation of Durhuus and
Jonsson.  Our definition of TQFT is category theoretic, as this seems
to be the most natural setting.  Section 3 shows that the vector
space associated with the $(d-1)$-sphere inherits the structure of a
commutative Frobenius algebra from the TQFT, and most importantly that
the decomposition of the TQFT into indecomposable theories is exactly the
decomposition of this Frobenius algebra into indecomposable
subalgebras.  It also classifies indecomposable Frobenius algebras in
terms of ordinary indecomposable algebras.
This, together with the parallel
theory for unitary TQFT's and $C^*$-Frobenius algebras, which is
developed along the way, gives a decomposition theory for TQFT's.
Section 4 shows that two-dimensional TQFT's are determined by
their Frobenius algebras, and gives a partial description of the
two-dimensional indecomposable TQFT's.  Section 5 ends with some
remarks.

I would like to thank Michael Artin, Scott Axelrod, John Baez, Lisa
Sawin, Is
Singer and Washington Taylor for help and conversations.  I would also
like
to thank John Barrett for pointing out to me the existence of
non-semisimple abelian Frobenius algebras.
\section{Axiomatic TQFT's and Direct Sums}

The cobordism category $\cob(d)$ in dimension $d$ has as objects closed,
oriented, $(d-1)$-dimensional smooth manifolds. A morphism with domain
$\Sigma_1$ and codomain $\Sigma_2$, called a cobordism from $\Sigma_1$
to $\Sigma_2$,  is up to diffeomorphism an oriented
$d$-dimensional
manifold with boundary $\Sigma_1^* \cup \Sigma_2$, where $\Sigma^*$ is
the manifold $\Sigma $ with the opposite orientation.  By `up to
diffeomorphism' we mean two morphisms $M$ and $M'$ are considered the
same if there is a boundary- and orientation-preserving diffeomorphism
between them.  More precisely, a morphism should be a $d$-manifold $M$
together with a choice of some subset of the boundary components to be
the domain, and a choice of ordering of the components of the domain
and codomain.  We will avoid this red tape and always make the domain,
codomain and order  clear from context.

Composition is by gluing: If $M_1:\Sigma_1 \to \Sigma_2$ and
$M_2:\Sigma_2 \to \Sigma_3$, then $M_2 M_1:\Sigma_1 \to \Sigma_3$ is
the manifold formed by identifying points in $M_1$ and $M_2$ on the
shared boundary $\Sigma_2$.  The identity $\bold{1}_\Sigma$ for each object
$\Sigma$ is
the cobordism $\Sigma \times I$.

Disjoint union gives a tensor product structure on $\cob(d)$.  That is,
we have a covariant functor from $\cob(d) \times \cob(d)$ to $\cob(d)$
which sends $\Sigma_1 \times \Sigma_2$  to $\Sigma_1 \cup \Sigma_2$,
and if $M:\Sigma_1 \to \Sigma_2$ and $M':\Sigma_1' \to \Sigma_2'$ then
$M  \times M'$ goes to $M \cup M': \Sigma_1 \cup \Sigma_1' \to
\Sigma_2 \cup \Sigma_2'$.
The empty $(d-1)$-manifold $\bold{\emptyset}$ is the
trivial object, and the empty $d$-manifold $\bold{1}_{\emptyset}$ is the
trivial morphism. That is,
they are the units for $\cup$ on objects and morphisms.  Furthermore
the cobordism $\frak{c}_{\Sigma_1,\Sigma_2}:\Sigma_1 \cup \Sigma_2 \to
\Sigma_2 \cup \Sigma_1$ given by the union of $\Sigma_1 \times I$
with $\Sigma_2 \times I$ with the boundary components ordered appropriately,
satisfies
$$\frak{c}_{\Sigma_1,\Sigma_2} \frak{c}_{\Sigma_2,\Sigma_1} =
\bold{1}_{\Sigma_1
  \cup \Sigma_2},$$
$$\frak{c}_{\Sigma_1 \cup \Sigma_2,\Gamma} = (\frak{c}_{\Sigma_1,\Gamma}
  \otimes \bold{1}_{\Sigma_2})(\bold{1}_{\Sigma_1} \otimes
  \frak{c}_{\Sigma_2,\Gamma})$$
and thus makes the cobordism category a symmetric monoidal or tensor
category \cite{MacLane71}.

$\mbox{Vect}(\FF)$ is also a tensor category.
The objects are finite-dimensional vector spaces over a field $\FF$,
and morphisms with domain $V$ and codomain $W$ are linear maps from $V$
to $W$.  Composition of morphisms is composition of linear maps, and
tensor product of objects and morphisms is tensor product of vector
spaces and linear maps.  Also $\bold{1}_V$ is the identity map on $V$, the
trivial object is $\FF$, and the trivial morphism is multiplication by
$1$ on $\FF$.

\begin{definition}\cite{Atiyah89,Atiyah90b}
A TQFT is a functor $\cal{Z}$ of tensor categories from
$\cob(d)$ to $\mbox{Vect}(\FF)$.  Two TQFT's are considered {\em
  equivalent\/} if there is a natural isomorphism between them.
\end{definition}

Of course, there is also a duality structure on cobordisms.  If
$M:\Sigma_1 \to \Sigma_2$, then we can consider $M^*$ as a cobordism
from $\Sigma_2 $ to $\Sigma_1$.  This kind of duality is analogous to
that in the category $\hil$ of Hilbert spaces and bounded linear
functionals, where each $f:H_1 \to H_2$ has an adjoint $f^*:H_2 \to
H_1$.  This motivates the following strengthening:

\begin{definition}
A {\em unitary TQFT\/} is a tensor functor from $\cob(d)$ to
$\hil$ such that $\cal{Z}(M^*)=
\cal{Z}(M)^*$.
\end{definition}

Durhuus and Jonsson \cite{DJ94} define the notion of the direct sum of
two TQFT's.  The direct sum of $\cal{Z}_1$ and $\cal{Z}_2$ is the
theory $\cal{Z}$ which associates to each connected $\Sigma$ the
vector space $\cal{Z}_1(\Sigma) \oplus \cal{Z}_2(\Sigma)$, associates
to each disconnected $\Sigma$ the tensor product of the vector spaces
associated to its components, associates to each connected $M$ the
linear map $\cal{Z}_1(M) \oplus \cal{Z}_2(M)$, interpreted in the
obvious way as an operator on the appropriate vector spaces, and
associates to each disconnected $M$ the tensor product of the values of
the components.  The reader may check that this is again a TQFT, and
that if $\cal{Z}_1$ and $\cal{Z}_2$ are unitary, so is the direct
sum.
\section{Frobenius Algebras and Decomposition of TQFT's}

Recall that a {\em Frobenius algebra\/} is a finite-dimensional
algebra $A$ over a field
$\FF$, together with a linear functional $\mu: A \to \FF$ such that
the bilinear pairing $(a,b)= \mu(ab)$ is nondegenerate.  If the pairing is
symmetric, for example if the algebra is commutative, we get what
Quinn \cite{Quinn95} calls an {\em ambialgebra\/}.  If $A$ is a
$C^*$-algebra and $\mu$ is a positive functional (i.e.,
$\mu(a^*a)>0$ for all nonzero $a \in A$) then we call $A$ a
$C^*$-Frobenius algebra.  The following
theorem is essentially due to Dijkgraaf \cite{Dijkgraaf89}.  Let $S$
be the $(d-1)$-sphere $S^{d-1}$.

\prop{pr:frobenius}
If $\cal{Z}$ is a $d$-dimensional TQFT, then $\cal{Z}$ gives
$\cal{Z}(S)$ the structure of a commutative Frobenius algebra,
and an action of this algebra on $\cal{Z}(\Sigma)$ for each
connected $(d-1)$-manifold $\Sigma$.  If $\cal{Z}$ is unitary,
$\cal{Z}(S)$ is  a $C^*$-Frobenius algebra  and the action
on $\cal{Z}(\Sigma)$ is a  $C^*$-representation.
\eprop

\begin{pf}

Let $A$ be $\cal{Z}(S)$. Multiplication is given by
the map $\cal{Z}(M_2^1):A
\tnsr A \to A$, where $M_2^1$ is the
$d$-ball with two $d$-balls removed.  The unit is the image of $1$
under the map $\cal{Z}(M_0^1):\FF \to A$, where $M_0^1$ is the
$d$-ball.  That it is a unit is immediate, associativity follows from
the fact that both sides of the associativity equation are $\cal{Z}$
of the $d$-ball with three $d$-balls removed.  Commutativity follows
from the fact that $M_2^1= M_2^1 \frak{c}_{S,S}$.

The map $\mu$ is $\cal{Z}(M_1^0):A \to \FF$, where $M_1^0$ is
$S^d$ with one $d$-ball removed.  The pairing is then $\cal{Z}(M_2^0)$,
where $M_2^0$ is the $d$-sphere with two balls removed.  To see that
it is nondegenerate, let $M_0^2$ be the connect sum of two $d$-balls,
and notice $(\bold{1}_S \cup M_2^0)(M_0^2 \cup \bold{1}_S)= \bold{1}_S$, so
that if
$\cal{Z}(M_0^2)(1)= \sum_i a_i \tnsr b_i \in A \tnsr
A$, then we have $\sum_i(x,a_i)b_i=x$, and thus the pairing
is nondegenerate.  See Figures \ref{fg:structure} and \ref{fg:axioms} for a
pictorial
presentation of the Frobenius algebra structure and axioms.

\fig
\begin{tabular*}{5.7in}{c@{\extracolsep{\fill}} c@{\extracolsep{\fill}}
    c@{\extracolsep{\fill}} c}
\pic{cup}{20}{-15}{0}{0}&
\pic{cap}{20}{-15}{0}{0}&
\pic{product}{30}{-20}{0}{0}&
\pic{pairing}{30}{-20}{0}{0}\\
 Unit & Mu & Product & Pairing
\end{tabular*}
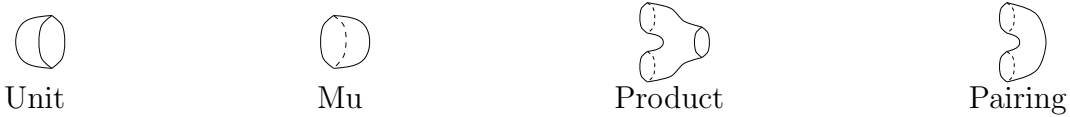
\efig{The structure of a Frobenius algebra}{fg:structure}
\fig
\begin{tabular*}{5.9in}{c@{\extracolsep{\fill}} c@{\extracolsep{\fill}}
    c}
$\pic{lassoc}{30}{-11}{0}{0} = \pic{massoc}{36}{-15}{0}{0}
= \pic{rassoc}{30}{-10}{0}{0}$ &
$\pic{lunit}{26}{-10}{0}{0}= \pic{runit}{18}{-10}{0}{0}$ &
$\pic{lnondeg}{30}{-10}{0}{0} =
\pic{rnondeg}{30}{-10}{0}{0}$\\
Associativity & Unit Axiom & Nondegeneracy
\end{tabular*}
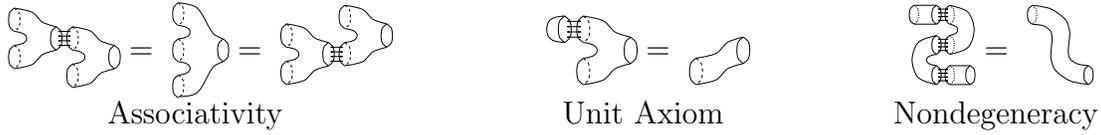
\efig{The  axioms of a Frobenius algebra}{fg:axioms}

For the action of $A$ on $\cal{Z}(\Sigma)$, let $\Sigma $ be
connected and let $M^\Sigma_{1,\Sigma}: S \cup \Sigma
\to \Sigma$ be $\Sigma \times I$ with a $d$-ball removed.  It is easy
to check that this is an algebra action.

Now suppose that $\cal{Z}$ is unitary.  The fact that the pairing is
nondegenerate means there is a conjugate-linear isomorphism $*$ from
$A$ to itself such that $(a^*,b)= \bracket{a,b}$, where
$\bracket{\cdot,\cdot}$ is the inner product on $A$.  Notice
that $(M^\Sigma_{1,\Sigma})^*= (1 \cup M^\Sigma_{1,\Sigma})(M_0^2 \cup 1)$, so
that
$\bracket{a  x,y}= \sum_i \bracket{a \tnsr x, a_i \tnsr b_iy} = \sum_i
\bracket{x,\bracket{a,a_i}b_iy} =  \sum_i
\bracket{x,(a^*,a_i)b_iy}=\bracket{x,a^* y}$.  Thus $A$ is
a $*$-algebra, and the representation onto
each $\cal{Z}(\Sigma)$ is a $*$-representation.
Since this representation on $\cal{Z}(S)$ is faithful, $A$ is a
$C^*$-algebra with the operator norm in this representation,
and the
representations on $\cal{Z}(\Sigma)$ will be $C^*$-algebra
representations.  Finally $\mu(a^*a)=\bracket{a,a} >0$, so $\mu$ is
positive.
\end{pf}
We say that $\cal{Z}$ is {\em based on \/} the commutative Frobenius
algebra $A$ if $\cal{Z}(S)$ is isomorphic to $A$ as a Frobenius
algebra.

\thm{th:decompose} Suppose
$\cal{Z}$ is based on a direct sum  $A=A_1 \oplus A_2$ of Frobenius algebras.
Then there exist
TQFT's $\cal{Z}_1$ and $\cal{Z}_2$, based on  $A_1$ and
$A_2$ respectively, such that $\cal{Z}= \cal{Z}_1\oplus\cal{Z}_2$.
Conversely, If $\cal{Z}$ decomposes as a direct sum of
theories, than the associated Frobenius algebra decomposes as a
corresponding direct
sum of Frobenius algebras.  Further, the same is true for unitary
TQFT's and direct sum of $C^*$-Frobenius algebras.
\ethm

\begin{pf}
Let $p_1$ and $p_2$ be the elements of $A$ which correspond
to the identities of $A_1$ and $A_2$ respectively, so that $p_ip_j=
\delta_{i,j}p_i$ and $p_1 + p_2 =1$.  Thus if we define
$\cal{Z}_i(\Sigma)$ for $\Sigma$ connected to be the range of the
action of $p_i$, we have that $\cal{Z}(\Sigma)= \cal{Z}_1(\Sigma)
\oplus \cal{Z}_2(\Sigma)$.  Likewise define the action of $p_i$ on
$\cal{Z}(\Sigma)$ for disconnected $\Sigma$ to be the tensor product
of its action on each connected component, and define
$\cal{Z}_i(\Sigma)$ to be its range.  Let $p_i$ act on
$\cal{Z}(\bold{\emptyset}) = \FF$ as $1$.  This defines $\cal{Z}$ on
$(d-1)$-manifolds in a manner satisfying the assumptions of a TQFT.

Now let $M:\Sigma \to \Gamma$ be a $d$-cobordism.  Define $M': S^{\cup
  k} \cup \Sigma \to \Gamma$ to be $M$ with a $d$-ball removed from
each of its $k$ components, and define $\cal{Z}_i(M)(x)=
\cal{Z}(M')(p_i^{\tnsr k}\tnsr x)$ for $x \in \cal{Z}(\Sigma)$.
Notice if we had removed more than one ball from some components and
put $p_i$'s in the appropriate tensor factors, we would have gotten
the same operator:  Removing two $d$-balls from the same component can
be regarded as removing one $d$-ball and gluing $M_2^1$ in, so that we
get the same operator as removing one ball and applying the result to
$p_i \cdot p_i=p_i$.  See Figure \ref{fg:redundant} for a pictorial
version of this argument.
\fig
$$ \pic{lredund}{60}{-30}{0}{0} \;=\;
\pic{mredund}{60}{-30}{0}{0} \;=\;
\pic{rredund}{50}{-25}{0}{0}$$
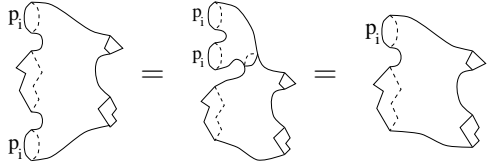
\efig{Removing two balls is the same as removing one}{fg:redundant}

We claim  that
$\cal{Z}_i(M)\cal{Z}_i(N)= \cal{Z}_i(MN)$.  Since this implies that
$\cal{Z}_i(M)\cal{Z}(M^\Sigma_{1,\Sigma})=\cal{Z}_i(M){\cal
      Z}_i(\bold{1}_\Sigma)= \cal{Z}_i(M)= {\cal
      Z}(M^\Gamma_{1,\Gamma})\cal{Z}_i(M)$, we have $\cal{Z}_i(M):
      \cal{Z}_i(\Sigma) \to \cal{Z}_i(\Gamma)$ and $\cal{Z}_i$ is a
      functor.  Since $\cal{Z}_i (M \cup N)= \cal{Z}_i(M) \tnsr
      \cal{Z}_i(N)$, it is a tensor functor and hence a TQFT.

To see this claim, notice that $\cal{Z}_i(M)\cal{Z}_i(N)(x)$ is
$\cal{Z}(K)(p_i^{\tnsr n}\tnsr x)$, where $K$ is $MN$ with some
positive number of $d$-balls removed from each component.   We have
already seen this is the same as $\cal{Z}_i(MN)$ (see Figure
\ref{fg:glue}).
\fig
$$ \pic{lglue}{60}{-30}{0}{0} \;=\;
\pic{mglue}{60}{-30}{0}{0} \;=\;
\pic{rglue}{50}{-25}{0}{0}$$
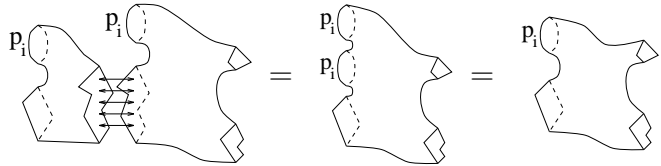
\efig{The action of $p_i$ intertwines composition}{fg:glue}

  Thus $\cal{Z}_i$
is  a TQFT.

Is $\cal{Z}= \cal{Z}_1 \oplus \cal{Z}_2$?  For $\Sigma$ connected we
have $\cal{Z}(\Sigma)= \cal{Z}_1(\Sigma) \oplus \cal{Z}_2(\Sigma)$,
and clearly if $\Sigma$ is not connected, then $\cal{Z}(\Sigma)$ is the
tensor product of $\cal{Z}$ of the connected components.  If $M$ is
connected, then $\cal{Z}_i(M)(x)= \cal{Z}(M')(p_i \tnsr x)$, so
$\cal{Z}_1(M)(x) + \cal{Z}_2(M)(x)=
\cal{Z}(M')((p_1 + p_2)\tnsr x)= \cal{Z}(M)(x)$.  Likewise, if $M$ is
disconnected then $\cal{Z}(M)$ is the tensor product of the values of
$\cal{Z}$ on the connected components.   Thus $\cal{Z}= \cal{Z}_1
\oplus \cal{Z}_2$.

The converse is easy.  If $\cal{Z}= \cal{Z}_1 \oplus \cal{Z}_2$, then
$\mu= \cal{Z}(M_1^0)= \cal{Z}_1(M_1^0) \oplus \cal{Z}_2(M_1^0)= \mu_1
\oplus \mu_2$.  Also $\bold{m} = \cal{Z}(M_2^1)= \cal{Z}_1(M_2^1) \oplus
\cal{Z}_2(M_2^1)= \bold{m}_1
\oplus \bold{m}_2$, where $\bold{m}$, $\bold{m}_1$, and $\bold{m}_2$ represent
the
products in $A$, $A_1$ and $A_2$ respectively.  Thus $A= A_1 \oplus A_2$.

If $\cal{Z}$ is unitary and its $C^*$-Frobenius algebra $A$
decomposes as $C^*$-Frobenius algebra into a direct sum $A_1 \oplus
A_2$, then the above
argument shows that $\cal{Z}=\cal{Z}_1 \oplus \cal{Z}_2$ as TQFT's:
We just need to show that $\cal{Z}_i$ is unitary.   But $p_i$ is a
self-adjoint projection, since the direct sum was as a $C^*$-algebra
direct sum,
so
\begin{align*}
\bracket{y,\cal{Z}_i(M) x} &= \bracket{y,\cal{Z}(M)(p_i x)}\\
&= \bracket{\cal{Z}(M)^*y,p_i x}\\
&= \bracket{p_i\cal{Z}(M^*) y,x}\\
&= \bracket{\cal{Z}(M^*)(p_i y),x}\\
&= \bracket{\cal{Z}_i(M^*)y,x}.
\end{align*}

Thus $\cal{Z}_i$ is unitary.  Conversely, if $\cal{Z}= \cal{Z}_1
\oplus \cal{Z}_2$ is a direct sum of unitary theories, then the
subspaces $A_1$ and $A_2$ of $A$ are orthogonal, so the $C^*$-norm is
the direct sum norm and the involution is the direct sum involution.
\end{pf}

Thus the direct sum decomposition of a TQFT corresponds exactly to the
direct sum decomposition of its associated commutative Frobenius
algebra.  Indecomposable commutative $C^*$-Frobenius algebras are easy
to classify:  For $\lambda \in \RR^+$, define the commutative $C^*$-Frobenius
algebra
$\bold{C}_\lambda$ to be the  $C^*$-algebra
$\CC$, with $\mu(x)=\lambda^{-1} x$.
It is clearly a simple $C^*$-algebra, and $\mu$ is
positive.  Since the only indecomposable commutative $C^*$-algebra is
one-dimensional, it is clear this exhausts all the possibilities.

\begin{corollary} Every unitary TQFT is a direct sum of unitary
  TQFT's, each  based on the Frobenius algebra $\bold{C}_\lambda$ for
  some $\lambda$.  In particular, every unitary TQFT is the direct sum
  of theories with one-dimensional $\cal{Z}(S^{d-1})$.
\end{corollary}

The story is a bit more complicated for arbitrary Frobenius algebras,
because there are so many commutative algebras.
Nevertheless, we can get a fairly complete description modulo this
issue.  Assume $\FF$ is algebraically closed.  For each $\lambda \in
\FF$ nonzero, let $\bold{S}_\lambda$ be the algebra $\FF$ with $\mu(x)=
\lambda^{-1} x$.  Also, let $A$ be a
commutative algebra spanned by the identity and at least one nilpotent, and
suppose the {\em socle\/}, the space of all $x \in A$ such that
$ax=0$ for all nilpotent $a \in A$, is  one-dimensional.
Let $\mu$ be  any linear functional
on $A$ which is nonzero on the socle.  Let $\bold{N}_{A,\mu}$ be this algebra
together
with this functional.

\prop{pr:Frobenius-decomp}  $\bold{S}_\lambda$ and $\bold{N}_{A,\mu}$
are indecomposable Frobenius algebras.  Further, every commutative
indecomposable Frobenius algebra is isomorphic to  one of
these, and these are nonisomorphic up to algebra isomorphism.
\eprop

\begin{pf}

Obviously $\bold{S}_\lambda$ is an indecomposable Frobenius algebra.
Now consider $\bold{N}_{A,\mu}$.  Notice for any finite-dimensional
algebra
there is a bound on the number of nilpotent elements  which can be
multiplied to get a nonzero product.  Thus for any $x \in A$, there must be
a $y\in A$ such that $xy$ is a nonzero element of the socle, for otherwise we
would be
able to write an arbitrarily long product of nilpotents times $x$ which
is nonzero.  For this $y$ we have $\mu(xy)$ is nonzero, and thus the
bilinear form is nodegenerate.  So $\bold{N}_{A,\mu}$ is a Frobenius
algebra.  It is clearly indecomposable, because $A$ is indecomposable.

To see that every indecomposable Frobenius algebra is one of the
above, first note that if
a Frobenius algebra decomposes into a direct sum as an algebra, it
decomposes in the same way as a Frobenius algebra, because the summands
are orthogonal subspaces with respect to the bilinear pairing, and
thus the bilinear pairing is nondegenerate when restricted to each.
So every indecomposable  Frobenius algebra is also
indecomposable as an algebra.   Thus it is spanned by the identity and
nilpotents.  If it has no nilpotents, then it is one-dimensional, and it is
clearly isomorphic to exactly one $\bold{S}_\lambda$.  If it has
nilpotents, then arguing as above it has a nonempty socle.   Every
element of the socle is orthogonal to all nilpotent elements, so the
socle is dual to the space spanned by the identity, and thus is
one-dimensional.  So the algebra is isomorphic to exactly one of the
algebras used to construct the $\bold{N}$'s.  Clearly $\mu$ must be
nonzero on the socle, or otherwise the socle would be orthogonal to
all of $A$.  Thus the algebra isomorphism extends to a Frobenius
algebra isomorphism to exactly one $\bold{N}_{A,\mu}$.
\end{pf}

We call a TQFT  {\em simple\/} if it is based on $\bold{S}_\lambda$, i.e., if
$\cal{Z}(S)$ is one-dimensional.
We call a TQFT based on $\bold{N}_{A,\mu}$ {\em nilpotent\/}.

\begin{corollary}
Every TQFT is
a direct sum of simple and nilpotent theories.
\end{corollary}
\section{TQFT's in Two Dimensions}

In two dimensions, the Frobenius algebra completely determines the
TQFT.  The following theorem is due to Dijkgraaf \cite{Dijkgraaf89}.

\thm{th:2D}  There is exactly one two-dimensional TQFT based on a given
commutative Frobenius algebra.  There is exactly one two-dimensional unitary
TQFT
based on a given commutative $C^*$-Frobenius algebra.
\ethm

\begin{pf}
  First let us construct the TQFT from a given commutative Frobenius
  algebra $A$.  Certainly the vector space associated to $n$ copies of
  $S^1$ is $A^{\tnsr n}$ and to the empty one-manifold is $\FF$. We
  must associate operators to two-manifolds.  Clearly $\cal{Z}(M^1_2)$
  is the product, $\cal{Z}(M^0_1)$ is  $\mu$, $\cal{Z}(M^1_0)$ is
  the identity, and $\cal{Z}(M_1^2)$ is the dual map to the product
  under the pairing.  If $M$ is any
  two-dimensional cobordism, use a Morse function to write it as a
  product of cobordisms, each of which is a union of some number of
  copies of $\bold{1}_S$ with one of $M^0_1$, $M^1_0$, $M^2_1$,
  $M^1_2$. The
  value of $\cal{Z}$ on such a cobordism is determined by the
  requirement of tensor functoriality and the values of these four
  cobordisms.  To check it does not depend on the Morse function,
  recall by Cerf theory \cite{Cerf70} that any change of Morse function has the
  effect of a sequence of the following moves, illustrated in Figure
\ref{fg:Morse}:
\begin{romanlist}
\item $M^m_n \bold{1}_{S^{\cup n}} = M^m_n = \bold{1}_{S^{\cup m}} M^m_n$
\item $(M^0_1 \cup \bold{1}_S)(M^2_1)  = (\bold{1}_S \cup
  M^0_1)(M^2_1) = \bold{1}_S =
  M^1_2(M^1_0 \cup \bold{1}_S) =
  M^1_2(\bold{1}_S \cup M^1_0)$
\item $(M^1_2 \cup \bold{1}_S)(\bold{1}_S \cup M^2_1) = M^2_1
  M^1_2 = (\bold{1}_S \cup
  M^1_2) (M^2_1 \cup \bold{1}_S).$
\end{romanlist}
\fig
\begin{tabular*}{5.9in}{c@{\extracolsep{\fill}} c@{\extracolsep{\fill}}
    c}
$\pic{lmorse1}{40}{-20}{0}{0}=\pic{rmorse1}{30}{-15}{0}{0}$
& $\pic{lmorse2}{40}{-20}{0}{0}=
\pic{mmorse2}{40}{-20}{0}{0} =
\pic{rmorse2}{40}{-20}{0}{0}$ &
$\pic{lmorse3}{40}{-20}{0}{0} =
\pic{rmorse3}{40}{-20}{0}{0}$ \\
(i) & (ii) & (iii)
\end{tabular*}
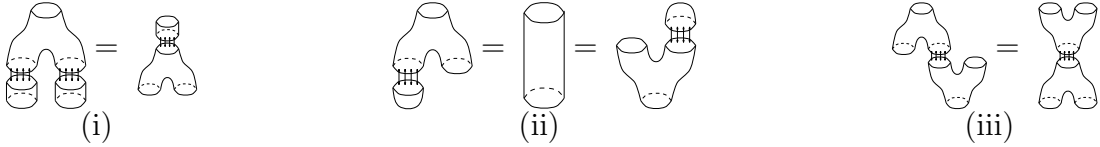
\efig{Moves of Cerf theory}{fg:Morse}

Our construction of $\cal{Z}(M)$ is clearly invariant
under $(i).$ Move $(ii)$ is just the statement that $1 \in A$ is the
identity, and that $\mu$ is dual to the identity.  For $(iii)$ notice
$\cal{Z}(M^2_1)(x)=\sum_i xa_i \tnsr b_i,$ where $a_i$ is a basis
of $A$ and $b_i$ is its dual basis.  Then the left-hand side of
$(iii)$ on $x \tnsr y$ is $\sum xy a_i \tnsr b_i = \cal{Z}(M^2_1)
\cal{Z}(M_1^2)(x \tnsr y)$, the right-hand side.  Likewise $\cal{Z}(M^2_1)(x) =
\sum a_i
\tnsr x b_i,$ showing the other equality.  Thus $\cal{Z}(M)$ is
well-defined.

A Morse function on $M_1$ and $M_2$ gives a Morse function on $M_1 M_2,$
so $\cal{Z}(M_1M_2)=\cal{Z}(M_1)\cal{Z}(M_2).$ Clearly
$\cal{Z}(\bold{1}_{S^{\cup n}}) = \bold{1}_{A^{\tnsr n}},$ so $\cal{Z}$ is a
functor.  It is easy to check that it takes union of one-manifolds and
cobordisms to tensor product of vector-spaces and operators
respectively, the empty one- and two-manifold get sent to $\FF$ and
$1$ respectively, and that permutation of components
corresponds to permutation of tensor factors.  Thus $\cal{Z}$ is a
TQFT.

If $\cal{Z}'$ is another TQFT based on the same Frobenius algebra $A$,
the identification of the Frobenius algebras gives a linear isomorphism
between the vector spaces of $\cal{Z}$ and $\cal{Z}'$.  What's more,
this map intertwines $\cal{Z}(M^m_n)$ with $\cal{Z}'(M^m_n)$ with
$m, n = 0,1,2$ as above, since their value is determined by the Frobenius
algebra.  The fact that $\cal{Z}$ and $\cal{Z}'$ are tensor functors
then ensures that the isomorphism intertwines $\cal{Z}(M)$ with
$\cal{Z}'(M).$ Since it obviously intertwines the tensor product
structure, this is a natural isomorphism.

If $A$ is a $C^*$-Frobenius algebra, then $\bracket{a,b}= \mu(a^*b)$
defines a positive-definite  inner product on $A$, and hence on every
vector space associated with the theory.  One need only check that
 $\cal{Z}(M_1^0)=\cal{Z}(M_0^1)^*$ and $\cal{Z}(M_1^2) = \cal{Z}(M_2^1)^*$.
\end{pf}

Thus there is a two-dimensional indecomposable TQFT associated to each
$\bold{S}_\lambda$ and each $\bold{N}_{A,\mu},$ and an indecomposable
unitary TQFT to each $\bold{C}_\lambda$,  and
every two-dimensional TQFT is a direct sum of these.

It is worth noting what these theories look like.  Following Durhuus
and Jonsson \cite{DJ94}, we check that in the theory $\cal{Z}_\lambda$
associated to $\bold{S}_\lambda,$ every vector space can be associated to
$\FF$ in such a way that $\cal{Z}_\lambda (M) = \lambda^{-\chi
  (M)/2},$ where $\chi (M)$ is the Euler number.

For the nilpotent case, we will of course not be able to give a simple
description of the whole TQFT $\cal{Z}_{A,\mu},$ but we can describe a
surprisingly large amount.

We can write  a chain of ideals $A=N_n \superset N_{n-1} \superset
\cdots \superset N_1,$  where $N_1$ is the socle, and each $N_k$ is
the preimage in $A$ of the socle of $A/N_{k-1}.$  Choose a basis for $a_i$ of
$A$ which restricts to a basis of each $N_k,$ and let $b_i$ be its
dual basis.  We claim $a_i b_i =s,$ where $s$ is the element of the
socle with $\mu(s)=1.$  To see this, notice if $a_i \in N_k$ and $y$
is nilpotent, then $ya_i \in N_{k-1},$ and hence can be written as a
linear combination of $a_j$ for $j \neq i.$   Thus $ya_i b_i =0,$ and
$a_i b_i$ is in the socle.  Since $\mu(a_i b_i)=1$, we have $a_ib_i=s$.

Now let $M=M^1_2 M^2_1,$ a twice-punctured torus.  Then
$$\cal{Z}(M)(x)=\sum_i xa_i b_i = \dim(A)xs=\dim(A)f(x)s$$
where $f(x)$ is the functional on $A$ which is 1 on the identity and
zero on the nilpotents (the unique homomorphism to $\FF$).  From this
one easily concludes that the torus with $n$ incoming punctures and
$m$ outgoing punctures is sent to the operator
$$x_i \tnsr \cdots \tnsr x_n \mapsto f(x_1) \cdots f(x_n)
\dim(A)s^{\tnsr m}$$
and any manifold of genus more than one gets sent to the operator $0$ on
the appropriate space.

\section{Remarks}

\begin{itemize}
\item It would be nice to find an action or state sum definition of
  the two-dimensional nilpotent TQFT's $\cal{Z}_{A,\mu}.$ There is no obvious
  impediment to this, but the surprising behavior of this TQFT would
  appear to make it difficult.
\item It has been asked (e.g. by \cite{DJ94}) whether two TQFT's which
  agree on all closed $d$-manifolds are naturally isomorphic.  The
  answer is no, even in two dimensions, if we do not restrict to
  unitary theories.  For any $\cal{Z}_{A,\mu}$ the sphere gets
  sent to $\mu(1),$ the torus gets sent to $\dim(A),$ and all others
  get sent to zero.  Clearly many different $A$ and $\mu$ give the
  same values for these, and since they correspond to nonisomorphic
  Frobenius algebras, they correspond to inequivalent TQFT's.  The
  values on closed $d$-manifolds does determine the TQFT for
  2-dimensional unitary and semisimple theories.
\item It is clear that the proofs of Proposition \ref{pr:frobenius}
  and Theorem \ref{th:2D} really only involve the category of
  cobordisms.  Thus it would be natural to express them as corollaries
  to purely topological theorems about this category. Specifically, we
  could define the notion of a Frobenius object in a tensor category,
  and a Frobenius action of one object on another.  Then Proposition
  \ref{pr:frobenius} follows from the statement that $S$ is a Frobenius object
  in $\cob(d)$ with a Frobenius action on each connected $\sum,$ and
  Theorem \ref{th:2D} from  the statement that $\cob(2)$ is the free
  tensor category generated by one Frobenius object.
\item If $\cal{Z}$ is a $d$ dimensional TQFT, and $X$ is an
  $r$-manifold for $r<d$, then $\cal{Z}$ and $X$ together naturally
  give a $(d-r)$-dimensional TQFT, which assigns to each
  $(d-r-1)$-manifold $\Sigma$ the vector space $\cal{Z}(\Sigma \times
  X)$, and to each $(d-r)$-cobordism $M$ the operator $\cal{Z}(M
  \times X)$.  In particular for each $(d-2)$-manifold we get a
  two-dimensional TQFT, which we can classify as in the previous
  section.  This classification should give important information
  about the TQFT in a simple format.  For example, if we take the
  Chern-Simons TQFT and $X=S^1$, we get a sum of simple TQFT's with
  $\lambda_i= S_{0,i}^{-1}$.
\end{itemize}

\bibliographystyle{alpha}

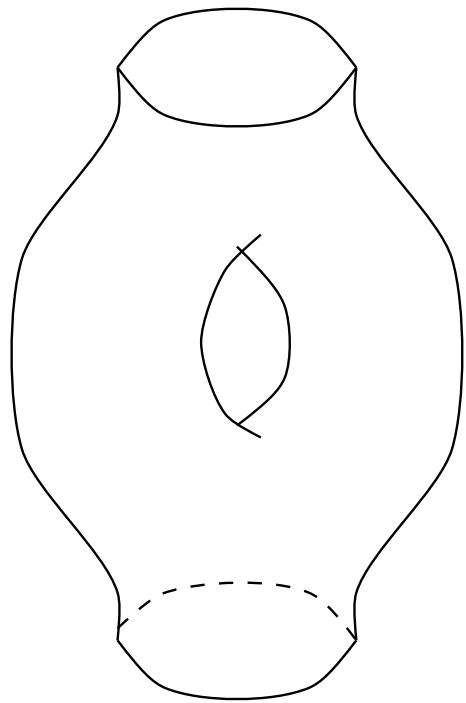
\end{document}